\documentclass{amsart}
\usepackage{amsmath} 
\usepackage{amssymb}
\usepackage{mathrsfs}

\newtheorem{theorem}{Theorem}[section] 
\newtheorem{claim}[theorem]{Claim}

\newtheorem{observation}[theorem]{Observation}

\theoremstyle{definition}
\newtheorem{definition}[theorem]{Definition}

\newtheorem{discussion}[theorem]{Discussion}
\newtheorem{convention}[theorem]{Convention}

\newtheorem{question}[theorem]{Question}

\theoremstyle{remark}
\newtheorem{remark}[theorem]{Remark}
\newtheorem{notation}[theorem]{Notation}
\newtheorem{conclusion}[theorem]{Conclusion}

\newcommand{\rest}{{\restriction}}
 
\newcommand{\odd}{{\rm odd}} 

\newcommand{\arity}{{\rm arity}}

\newcommand{\LST}{{\rm LST}}

\newcommand{\wilog}{{\rm without loss of generality}}

\newcommand{\then}{{\underline{then}}}
\newcommand{\when}{{\underline{when}}}
\newcommand{\Then}{{\underline{Then}}}

\newcommand{\Iff}{{\underline{iff}}}
\newcommand{\mn}{{\medskip\noindent}}
\newcommand{\sn}{{\smallskip\noindent}}

\newcommand{\cA}{{\mathcal A}}

\newcommand{\gB}{{\mathfrak B}}
\newcommand{\gh}{{\mathfrak h}}

\newcommand{\cH}{{\mathcal H}}

\newcommand{\bbL}{{\mathbb L}}
\newcommand{\cM}{{\mathcal M}}

\newcommand{\cP}{{\mathcal P}}

\newcommand{\gk}{{\mathfrak k}}
\newcommand{\gK}{{\mathfrak K}}

\newcommand{\cU}{{\mathcal U}}
\newcommand{\cX}{{\mathcal X}}

\newcommand{\cf}{{\rm cf}}

\newcount\skewfactor
\def\mathunderaccent#1#2 {\let\theaccent#1\skewfactor#2
\mathpalette\putaccentunder}
\def\putaccentunder#1#2{\oalign{$#1#2$\crcr\hidewidth
\vbox to.2ex{\hbox{$#1\skew\skewfactor\theaccent{}$}\vss}\hidewidth}}

\newenvironment{PROOF}[2][\proofname.]
   {\begin{proof}[#1]\renewcommand{\qedsymbol}{\textsquare$_{\rm #2}$}}
   {\end{proof}}

\begin{document}

\title {Maximal Failures of Sequence Locality in A.E.C.}
\author {Saharon Shelah}
\address{Einstein Institute of Mathematics\\
Edmond J. Safra Campus, Givat Ram\\
The Hebrew University of Jerusalem\\
Jerusalem, 91904, Israel\\
 and \\
 Department of Mathematics\\
 Hill Center - Busch Campus \\ 
 Rutgers, The State University of New Jersey \\
 110 Frelinghuysen Road \\
 Piscataway, NJ 08854-8019 USA}
\email{shelah@math.huji.ac.il}
\urladdr{http://shelah.logic.at}
\thanks{The author thanks Alice Leonhardt for the beautiful typing.
  First typed October 20, 2007. Paper 932}






\date{December 1, 2017}

\begin{abstract}
We are interested in examples of a.e.c. with amalgamation
having some (extreme) behaviour concerning types.  Note we deal with
${\frak k}$ being
sequence-local, i.e. local for increasing chains of length a regular
cardinal (for types, equality of all restrictions imply equality).
For any cardinal $\theta \ge \aleph_0$ we construct an
a.e.c. with amalgamation ${\frak k}$ with $\LST({\gk}) = 
\theta,|\tau_{\frak K}| = \theta$ such that
$\{\kappa:\kappa$ is a regular cardinal and ${\frak K}$ is not
$(2^\kappa,\kappa)$-sequence-local$\}$ is maximal.  In fact we have a direct
characterization of this class of cardinals: the regular $\kappa$ such
that there is no uniform $\kappa^+$-complete ultrafilter (on any
$\lambda > \kappa$).   We also
prove a similar result to ``$(2^\kappa,\kappa)$-compact for types".
\end{abstract}

\maketitle
\numberwithin{equation}{section}
\setcounter{section}{-1}
\newpage

\section {Introduction}

Recall a.e.c. (abstract elementary classes); were introduced in
\cite{Sh:88}; and their (orbital) types defined in \cite{Sh:300}, see on them
\cite{Sh:h}, \cite{Bal0x}.  It has seemed to me obvious that even with
${\frak k}$ having amalgamation, those types in general lack the good
properties of the classical types in model theory.
E.g. ``$(\lambda,\kappa)$-sequence-locality where

\begin{definition}
\label{z1}  
1) We say that an a.e.c. ${\gk}$
is a $(\lambda,\kappa)$-sequence-local (for types)
\when \, $\kappa$ is regular and for every
$\le_{\gk}$-increasing continuous sequence $\langle M_i:i \le
\kappa \rangle$ of models of cardinality $\lambda$ and $p,q 
\in \mathscr{S}(M_\kappa)$ we have $(\forall i < \kappa)(p
\restriction M_i = q \restriction M_i) \Rightarrow p = q$.  We omit
$\lambda$ when we omit ``$\|M_i\| = \lambda$".

\noindent
2) We say an a.e.c. ${\gk}$ is $(\lambda,\kappa)$-local
when: $\kappa \ge \text{ LST}({\frak k})$ and if $M \in 
{\gk}_\lambda$ and $p_1,p_2 \in \mathscr{S}(M)$ and $N \le_{\gk} M
\wedge \|N\| \le \kappa \Rightarrow p_1 \rest N = p_2 \rest N$
\then \, $p_1 = p_2$.

\noindent
3) We may replace $\lambda$ by $\le \lambda,< \lambda,[\mu,\lambda]$
with the obvious meaning (and allow $\lambda$ to be infinity).
\smallskip

Of course, being sure is not a substitute for a proof, some examples
were provided by Baldwin-Shelah \cite[\S2]{BlSh:862}.  There we give an
example of the failure of $(\lambda,\kappa)$-sequence-locality for
${\gk}$-types in ZFC for some $\lambda,\kappa$, actually $\kappa =
\aleph_0$.  This was done by translating our 
problems to abelian group problems.  While those
problems seem reasonable by themselves they may hide our real problem.

Here in \S1 we get ${\gk}$, an a.e.c. with amalgamation with the class
$\{\kappa:(< \infty,\kappa)$-sequence-localness fail for $\gk\}$ being 
maximal; what seems to me a major missing point up to it, see Theorem
\ref{e3}.  Also we
deal with ``compactness of types" getting unsatisfactory results -
classes without amalgamation; in \cite{BlSh:862} this was done
only in some universes of set theory but with amalgamation; see \S2.

We relay on \cite{BlSh:862} to get that ${\gk}$ has the JEP and
amalgamation. 
\end{definition}

\begin{question}
\label{z2}  Can $\{\kappa:{\gk}$ is $(< \infty,\kappa)$-local$\}$ be
``wild"?  E.g. can it be all odd regular alephs? etc?

Note that for this the present translation theorem of \cite{BlSh:862}
is not suitable.

In \S2 we deal with sequence-compactness of types.

We thank Will Boney for a correction.
\end{question}
\newpage

\section {An a.e.c. with maximal failure of being local} 

\begin{claim}
\label{e1}  
Assume
\mn
\begin{enumerate}
\item[$\circledast_1$]   $(a) \quad \kappa = \text{\rm cf}(\kappa) >
\theta \ge \aleph_0$ \underline{or just} $\kappa = \text{\rm cf}(\kappa) \ge
\aleph_0,\theta \ge \aleph_0$
\sn
\item[${{}}$]   $(b) \quad$ there is no uniform $\theta^+$-complete
ultra-filter $D$ on $\kappa$
\sn
\item[${{}}$]   $(c) \quad \tau_\theta$ is a vocabulary of
cardinality $\theta$ consisting of $\theta \, n$-place predicates

\hskip25pt  with each $n$ (and no more say $\{R_{\gamma,n}:\gamma <
\theta,n < \omega\},n = \arity(R_{\gamma,n}))$.
\end{enumerate}
\mn
\Then
\mn
\begin{enumerate}
\item[$\boxplus$]  there are $I_\alpha,M_{\ell,\alpha},
\pi_{\ell,\alpha}$ (for $\ell=1,2$ and $\alpha \le \kappa$), 
$g_\alpha$ (for $\alpha < \kappa$) satisfying:
\begin{enumerate}
\item[$(a)$]   $I_\alpha$, a set of cardinality $\theta^\kappa$,
 is $\subseteq$-increasing continuous with $\alpha$
\sn
\item[$(b)$]   $M_{\ell,\alpha}$,  a $\tau_\theta$-model
of cardinality $\le \theta^\kappa$, is increasing continuous with $\alpha$ 
\sn
\item[$(c)$]  $\pi_{\ell,\alpha}$ is a function from
$M_{\ell,\alpha}$ onto $I_\alpha$, increasing continuous with $\alpha$
\sn
\item[$(d)$]  $|\pi_{\ell,\alpha}^{-1}\{t\}| \le \theta^{\aleph_0}$ for $t
\in I_\alpha,\alpha \le \kappa$ and $\ell = 1,2$
\sn
\item[$(e)$]  if $t \in I_{\alpha +1} \backslash I_\alpha$ then
$\pi^{-1}_{\ell,\alpha}\{t\} \subseteq M_{\ell,\alpha +1} \backslash
M_{\ell,\alpha}$ 
\sn
\item[$(f)$]  for $\alpha < \kappa,g_\alpha$ is an isomorphism
from $M_{1,\alpha}$ onto $M_{2,\alpha}$ respecting
$(\pi_{1,2},\pi_{2,\kappa})$ which means $a \in M_{1,\alpha}
\Rightarrow \pi_{1,\alpha}(a) = \pi_{2,\alpha}(g_\alpha(a))$
\sn
\item[$(g)$]   for $\alpha = \kappa$ there is no isomorphism from
$M_{1,\alpha}$ onto $M_{2,\alpha}$ respecting
$(\pi_{1,\alpha},\pi_{2,\alpha})$.
\end{enumerate}
\end{enumerate}
\end{claim}

\begin{PROOF}{\ref{e1}}
Follows from \ref{e2} which is just a fuller 
version adding to $\tau_\theta$ unary function $F_c$ for $c \in G$;
this is just a notational change when $\theta^{\aleph_0} = \theta$.
Otherwise see $(*)_{12}$ of the proof of \ref{e2}; 
anyhow we shall use \ref{e2}.
\end{PROOF}

\begin{claim}
\label{e2}  
Assuming $\circledast_1$ of \ref{e1} we have:
\mn
\begin{enumerate}
\item[$\boxplus$]   there are $I_\alpha,M_{\ell,\alpha},
\pi_{\ell,\alpha}$ (for $\ell=1,2,\alpha \le \kappa$) and 
$g_\alpha$ (for $\alpha < \kappa$) and $G$ such that:
\begin{enumerate}
\item[$(a)$]   $G$ is an additive (so abelian) group of
cardinality $\theta^{\aleph_0}$
\sn
\item[$(b)$]   $I_\alpha$ is a set, increasing continuous with
$\alpha,|I_\alpha| = \theta^\kappa$
\sn
\item[$(c)$]  $M_{\ell,\alpha}$ is a $\tau^+_\theta$-model,
increasing continuous with $\alpha$, of cardinality $\theta^\kappa$
where $\tau^+_\theta = \tau_\theta \cup \{F_c:c \in G\},F_c$ a unary
function symbol, $\tau_\theta$ is from $\circledast_1(c)$ of \ref{e1}
\sn
\item[$(d)$]   $\pi_{\ell,\alpha}$ is a function from
$M_{\ell,\alpha}$ onto $I_\alpha$ increasing continuous with $\alpha$
\sn
\item[$(e)$]  $F^{M_{\ell,\alpha}}_c(c \in G)$ is a permutation of
$M_{\ell,\alpha}$, increasing continuous with $\alpha$
\sn
\item[$(f)$]   $\pi_{\ell,\alpha}(a) =
\pi_{\ell,\alpha}(F^{M_{\ell,\alpha}}_c(a))$ 
\sn
\item[$(g)$] $F^{M_{\ell,\alpha}}_{c_1}(F^{M_{\ell,\alpha}}_{c_2}(a)) =
F^{M_{\ell,\alpha}}_{c_1 + c_1}(a)$
\sn
\item[$(h)$]   $\pi_{\ell,\alpha}(a) = \pi_{\ell,\alpha}(b)
\Leftrightarrow \bigvee\limits_{c \in G} F^{M_{\ell,\alpha}}_c(a)=b$
\sn
\item[$(i)$]   for $\alpha < \kappa,g_\alpha$ is an isomorphism
from $M_{1,\alpha}$ onto $M_{2,\alpha}$ which respects
$(\pi_{1,\alpha},\pi_{2,\alpha})$ which means
$a \in M_{1,\alpha} \Rightarrow \pi_{1,\alpha}(a) =
\pi_{2,\alpha}(f_\alpha(a))$
\sn
\item[$(j)$]   there is no isomorphism from $M_{1,\kappa}
\restriction \tau_\theta$ onto $M_{2,\kappa} \restriction \tau_\theta$ 
respecting $(\pi_{1,\kappa},\pi_{2,\kappa})$.
\end{enumerate}
\end{enumerate}
\end{claim}

\begin{PROOF}{\ref{e2}}  
Let
\mn
\begin{enumerate}
\item[$(*)_0$]   $\sigma = \theta^{\aleph_0}$ so $\sigma = \sigma^{\aleph_0}$
\sn
\item[$(*)_1$]  
\begin{enumerate}
\item[(a)]  let $G = ([\sigma]^{< \aleph_0},\Delta)$, i.e., the
family of finite subsets of $\sigma$ with the 
operation of symmetric difference.  This is an abelian group
satisfying $\forall x (x + x=0)$
\sn
\item[(b)]   let $\langle a_{f,\alpha,u}:f \in
{}^\kappa \sigma,\alpha < \kappa,u \in G\rangle$ be a sequence without
repetitions
\sn
\item[(c)]   for $\beta \le \kappa$ let $A_\beta =
\{a_{f,\alpha,u}:f \in {}^\kappa \sigma,\alpha < 1 + \beta$ and $u \in G\}$
\sn
\item[(d)]   for $\beta \le \kappa$
let $I_\beta = ({}^\kappa \sigma) \times (1 + \beta)$
\sn
\item[(e)]  let $\pi_\beta(a_{f,\alpha,u}) = (f,\alpha)$ when
$\alpha< 1 + \beta \le \kappa$
\sn
\item[(f)]  for each $\beta < \kappa$ we define a permutation
$g_\beta$ (of order 2) of $A_\beta$ by $g_\beta(a_{f,\alpha,u}) = 
a_{f,\alpha,u +_G\{f(\beta)\}}$ hence 
$a \in A_\beta \Rightarrow \pi_\beta(g_\beta(a)) = \pi_\beta(a)$.
\end{enumerate}
\end{enumerate}
\mn
Note that
\mn
\begin{enumerate}
\item[$(*)_2$] 
\begin{enumerate}
\item[(a)]  $|G| = \sigma$
\sn
\item[(b)]  $\langle A_\beta:\beta \le \kappa\rangle$
is a $\subseteq$-increasing continuous, each $A_\beta$ a set of cardinality 
$\sigma^\kappa = \theta^\kappa$
\sn
\item[(c)]  $\langle I_\beta:\beta \le \kappa\rangle$
is $\subseteq$-increasing continuous, each $I_\beta$ of cardinality
$\sigma^\kappa = \theta^\kappa$
\sn
\item[(d)]  $\pi_\beta$ is a mapping from $A_\beta$ onto $I_\beta$
\sn
\item[(e)]   if $t \in I_\alpha \subseteq I_\beta$
then $\pi^{-1}_\beta\{t\} = \pi^{-1}_\alpha\{t\}$ has cardinality $|G|
= \sigma$
\sn
\item[(f)]   if $t \in I_{\alpha +1} \backslash
I_\alpha$ then $\pi^{-1}_{\alpha +1}\{t\} \subseteq A_{\alpha +1}
\backslash A_\alpha$
\sn
\item[(g)]  if $\alpha \le \beta \le \kappa$ then $g_\beta$ maps
  $A_\alpha$ onto itself and $g_\beta \circ g_\beta$ is the identity.
\end{enumerate}
\end{enumerate}
\mn
For each $n < \omega$ and $\beta \le \kappa$ 
we define equivalence relations $E'_{n,\beta},E_{n,\beta}$ on ${}^n(A_\beta)$:
\mn
\begin{enumerate}
\item[$(*)_3$]   $\bar a E'_{n,\beta} \bar b$ iff $\pi_\beta(\bar a) =
\pi_\beta(\bar b)$ where of course $\pi_\beta(\langle a_\ell:\ell <
n\rangle) = \langle \pi_\beta(a_\ell):\ell < n \rangle$
\sn
\item[$(*)_4$]   $\bar a E_{n,\beta} \bar b$ \Iff \, $\bar a
E'_{n,\beta} \bar b$ and there are $k < \omega$ and $\bar a_0,\dotsc,\bar
a_k$ such that
\sn
\begin{enumerate}
\item[(i)]   $\bar a_\ell \in {}^n(A_\beta)$
\sn
\item[(ii)]  $\bar a = \bar a_0$
\sn
\item[(iii)]  $\bar b = \bar a_k$
\sn
\item[(iv)]  for each $\ell < k$ for some $\alpha_1,\alpha_2 <
\kappa$ we have $g^{-1}_{\alpha_2}(g_{\alpha_1}(\bar a_\ell))$ is well
defined and equal to $\bar a_{\ell +1}$ \underline{or}
$g_{\alpha_2}(g^{-1}_{\alpha_1}(\bar a_\ell))$ is well defined and
equal to $\bar a_{\ell +1}$.
\end{enumerate}
\end{enumerate}
\mn
Note: 
\mn
\begin{enumerate}
\item[$(*)_{4.1}$]
\begin{enumerate}
\item[(a)]   the two possibilities in $(*)_4(iv)$ are
one as $g^{-1}_\alpha = g_\alpha$ so the first one is a special case
of the second;
\sn
\item[(b)]   $g_\alpha$ \underline{does not preserve} 
$\bar a/E_{n,\beta}!$, in fact,
$a,g_\alpha(a)$ are never $E_{n,\beta}$ equivalent; 
\sn
\item[(c)]  clearly they are well defined
\underline{iff} $(\forall \ell \le k)[\bar a_\ell \in
{}^n(A_{\text{min}\{\alpha_1,\alpha_2\}})]$ because if $\alpha \le \beta$
then $g_\beta$ maps $A_\beta$ onto itself because
$g(a_{f_1,\alpha,u_2}) = a_{f_1,\alpha,u_2} \Rightarrow |u_1| + 1 =
|u_2| \mod 2$
\sn
\item[(d)]  if $\alpha \le \beta,a \in A_\alpha$, then $g_\beta$ maps
$a/E_{n,\beta}$ onto itself
\sn
\item[(e)]  if $\alpha,\beta \le \kappa$, \then \, $g_\alpha,g_\beta$
  commute (on the intersection of their domains, $A_{\min\{\alpha,\beta\}}$.
\end{enumerate}
\end{enumerate}
\mn
[Why?  E.g. for clause (b) note clause (d).]

Note
\mn
\begin{enumerate}
\item[$(*)_5$] 
\begin{enumerate}
\item[(a)]  $E'_{n,\beta},E_{n,\beta}$ are indeed
equivalence relations on ${}^n(A_\beta)$
\sn
\item[(b)]  $ E_{\beta,n}$ refine $E'_{\beta,n}$
\sn
\item[(c)]   if $n < \omega,\bar a \in {}^n(A_\beta)$
then $\bar a /E'_{n,\beta}$ has at most $\sigma$ members (really
exactly two but we shall use only its having $\le 2^\sigma$ members)
\sn
\item[(d)]   if $\alpha < \beta \le \kappa$ \then \,
$E'_{n,\beta} \restriction {}^n(A_\alpha) = E'_{n,\alpha}$ and
$E_{n,\beta} \restriction {}^n(A_\alpha) = E_{n,\alpha}$ 
(read $(*)_4(iv)$ carefully!)
\sn
\item[(e)]   if $\alpha < \beta \le \kappa,\bar a \in
{}^n(A_\alpha)$ and $\bar b \in \bar a/E'_{n,\beta}$
\then \, $\bar b \in {}^n(A_\alpha)$
\sn
\item[(f)]   if $g_\alpha(\bar a_\ell) = \bar b_\ell$
for $\ell=1,2$ then: $\bar a_1 E'_{n,\beta} \bar a_2$ \Iff \, 
$\bar b_1 E'_{n,\beta} \bar b_2$.
\end{enumerate}
\end{enumerate}
\mn
Now we choose a vocabulary $\tau^*_\theta$ of cardinality $2^\sigma$
(but see $(*)_{12}$) and for $\alpha \le \kappa$ we choose a 
$\tau^*_\theta$-model $M_{1,\alpha}$ such that:
\mn
\begin{enumerate}
\item[$(*)_6$]
\begin{enumerate}
\item[(a)]  $M_{1,\alpha}$ increasing with $\alpha$ with universe $A_\alpha$
\sn
\item[(b)]   assume that $\bar a,\bar b$ are 
$E'_{n,\alpha}$-equivalent (so $\bar a,\bar b \in {}^n(A_\alpha)$ and 
$\pi_\alpha(\bar a) = \pi_\alpha(\bar b))$; \then \, $\bar a,\bar b$
realize the same quantifier free type 
in $M_{1,\alpha}$ iff $\bar a E_{n,\alpha} \bar b$
\sn
\item[(c)]  $\tau^*_\theta = \{E_n:n < \omega\} \cup
  \{F_c:c \in G\} \cup \{R_e:e \in {}^\sigma \sigma\} \cup \{R_{n,i}:i
 < \sigma,n < \omega\}$ where $E_n,R_e$ are two-place predicates, 
$F_c$ a unary function symbol, $R_{n,i}$ is $n$-place predicate
\sn
\item[(d)]   for every function $e \in {}^\sigma\sigma$

\hskip25pt $R^{M_{1,\alpha}}_e = \{(a_{f_1,\beta_1,u_1},
a_{f_2,\beta_2,u_2}) \in A_\alpha \times A_\alpha:f_1 = e \circ f_2$ and

\hskip150pt if $i < \sigma$ then $i \in u_1$

\hskip150pt iff $(|\{j \in u_2:e(j)=i\}|$ is odd)$\}$
\newline
recalling $f_\ell \in {}^\kappa \sigma$
\sn
\item[(e)]  $E^{M_{1,\alpha}}_n = E_{n,\alpha}$ and
  $F^{M_{1,\alpha}}_c = F_c$ is defined by $F_\ell:A_\alpha
  \rightarrow A_\alpha$ satisfies $F_c(a_{f,\alpha,u}) :=
  a_{f,\alpha,u+_G c}$
\sn
\item[(f)]   if $\alpha \le \beta_\ell < \kappa$ 
  for $\ell=1,2$ \then \, $g^{-1}_{\beta_2} g_{\beta_1} \rest 
A_\alpha$ is an automorphism of $M_{1,\alpha}$.
\end{enumerate}
\end{enumerate}
\mn
[Why is this possible?  First, we shall show that 
for each $\alpha < \kappa,g_\alpha$
maps $R^{M_{1,\alpha}}_e$ onto itself. 

Assume we are given a pair
$(a_{f_1,\beta_1,u_1},a_{f_2,\beta_2,u_2})$ from $A_\alpha \times
A_\alpha$ so $\beta_1,\beta_2 < 1 + \alpha$ and $f_1 = e \circ f_2$ so
\mn
\begin{enumerate}
\item[$(*)_{6.1}$]  $(a_{f_1,\beta_1,u_1},a_{f_2,\beta_2,u_2}) \in
R^{M_{1,\alpha}}_e$ iff $u_1 = \{e(j):j \in u_2$ and
$(\exists^{\text{odd}} \iota \in u_2)(e(\iota) = e(j))$.
\end{enumerate}
\mn
[Why?  Read $(*)_6$ carefully, in particular note that 
if $i \notin \{e(j):j \in u_2\}$ then $i \notin u_1$.]
\mn
\begin{enumerate}
\item[$(*)_{6.2}$]  $(g_\alpha(a_{f_1,\beta_1,u_1}),
g_\alpha(a_{f_2,\beta_2,u_2})) \in R^{M_{1,\alpha}}_e$ iff 
\newline
$(a_{f_1,\beta_1,u_1 + \{f_1(\alpha)\}},
a_{f_2,\beta_2,u_2+_G\{f_2(\alpha)\}}) \in R^{M_{1,\alpha}}_e$ iff
\newline
$u_1 +_G \{f_1(\alpha)\} = \{e(j):j \in u_2 +_G \{f_2(\alpha)\}$ and
$(\exists^{\text{odd}} \iota \in (u_2 +_G f_2(\alpha))(e(\iota) = e(j))\}$.
\end{enumerate}
\mn
[Why?  Inside $(*)_{6.2}$ the first ``iff" holds by the definition 
of $g_\alpha$, the second ``iff" holds as in $(*)_{6.1}$.]

But $f_1 = e \circ f_2$ hence 
\mn
\begin{enumerate}
\item[$(*)_{6.3}$]  $f_1(\alpha) = e(f_2(\alpha))$
\sn
\item[$(*)_{6.4}$]  letting $x = f_2(\alpha) < \sigma$ we have $u_1 =
  \{e(j):j \in u_2$ and $(\exists^{\odd} \iota \in u_2)(e(\iota) =
  e(j))\}$ \Iff \, $u_1 +_G\{e(x)\} = \{e(j):j \in u_2 +_G \{x\}$ and
  $\exists^\odd \iota \in (u_2 +_G \{x\})(e(\iota)=e(j))\}$.
\end{enumerate}
\mn
[Why?  Check by cases according to whether $x \in u_2$ and whether
$e(x) \in u_1$.  I.e. by ``$G$ is of order two" it suffices to prove
the ``only if" so assume the first equality in $(*)_{6.2}$.  If $e(x)
\notin u_1$, then just add $e(x)$ to both sides.  Similarly if $e(x)
\in u_1 \cap X \notin u_2$ and if $e(x) \in u_1 \wedge x \in u_2$.]

So together we get equivalence, so the ``first" holds.

Second, for defining the $R^{M_{1,\alpha}}_{n,i}$'s
\mn
\begin{enumerate}
\item[$(*)_{6.5}$]
\begin{enumerate}
\item[(a)]  for each $N$ let $E''_n$ be the following equivalence
  relation on ${}^n G:\bar u_1 E''_n \bar u_2$ \Iff \, for some $v \in
  G,|v|$ is even and $\bigwedge\limits_{\ell < n} u_{1,\ell} +_G v =
  u_{2,\ell}$
\sn
\item[(b)]   let $\langle \cX_{n,i}:i < \sigma\rangle$ list 
the $E''_n$-equivalence classes
\sn
\item[(c)]  let $R^{M_{1,\alpha}}_{n,i} = \{\bar a \in
  {}^n(A_\alpha)$: if $\bar a = \langle
  a_{f_\ell,\alpha_\ell,u_\ell}:\ell < n\rangle$, \then \, $\langle
  u_\ell:\ell < n\rangle \in \cX_{n,i}\}$.
\end{enumerate}
\end{enumerate}
\mn
This completed the choice of $M_{1,\alpha}$.  Third, $g_\alpha$ preserves
``$\bar a,\bar b$ are $E_{n,\alpha}$-equivalent", ``$\bar a,\bar b$
are $E'_{n,\alpha}$-equivalent" and their negations.  That is, $\bar
a,g_\alpha(\bar a)$ are not $E_{n,\alpha}$-equivalent, but as
$(\forall \beta)(g_\beta = g_\beta^{-1}),\bar a,\bar b$ being
$E_{n,\alpha}$-equivalent means that there is an even length pass from
$\bar a$ to $\bar b$, in the graph $\{(\bar c,g_\beta(\bar c)):\beta \in
[\gamma,\kappa)$ and $\bar c \in {}^n(A_\gamma)\}$ where $\gamma =
\text{ min}\{\gamma:\bar a,\bar b \in {}^n(A_\gamma)\}$.  

Fourth, no problem in the $M_{1,\alpha}$'s are increasing by
$(*)_5(d)$, just check that.

Fifth, $g_\alpha$ commutes with $F^{M_{1,\alpha}}_c$ for $c \in G$
because $G$ is an Abelian group.

Sixth, we should check clause $(*)_6(f)$.  Now $g^{-1}_{\beta_2}
g_{\beta_1} \rest A_\alpha = (g_{\beta_2} \rest A_\alpha)(g_{\beta_1}
\rest A_\alpha)$ by $(*)_2(g)$ and it has order 2 because $G$ is of
order 2 and it maps $E^{M_{1,\alpha}}_n$ to itself by the ``third",
commute with $F^{M_{1,\alpha}}_c$ by the fifth, maps
$R^{M_{1,\alpha}}_e$ to itself by the ``first".

Lastly, it maps $R^{M_{1,\alpha}}$ to itself by $(*)_{6.4}$.  So we
are done proving $(*)_6$.
\mn
\begin{enumerate}
\item[$(*)_7$]   for $\alpha < \kappa$ let $M_{2,\alpha}$ be the
$\tau^*_\theta$-model with universe $A_\alpha$ such that $g_\alpha$
is an isomorphism from $M_{1,\alpha}$ onto $M_{2,\alpha}$.
\end{enumerate}
\mn
Now we note
\mn
\begin{enumerate}
\item[$(*)_8$]   for $\alpha < \beta < \kappa,M_{2,\alpha} \subseteq
M_{2,\beta}$.
\end{enumerate}
\mn
[Why?  By the definitions of $M_{1,\gamma},g_\gamma,E'_{n,\gamma},
E_{n,\gamma}$, in particular, the ``first" and ``third",
in ``why $(*)_6$", fourth, i.e. $(*)_5(d)$ .]
\mn
\begin{enumerate}
\item[$(*)_9$]   let $M_{2,\kappa} := \cup\{M_{2,\alpha}:\alpha <
\kappa\}$, well defined by $(*)_8$
\sn
\item[$(*)_{10}$]  let $\pi_{\ell,\beta} = \pi_\beta$ for
$\ell=1,2$ and $\beta < \kappa$
\sn
\item[$(*)_{11}$]   except clause (j) the
demands in the conclusion of $\boxplus$ of \ref{e2} holds easily.
\end{enumerate}
\mn
[Why?  Just check.]
\mn
\begin{enumerate}
\item[$(*)_{12}$]   it is O.K. to use a vocabulary of cardinality
$2^\sigma = 2^{\theta^{\aleph_0}}$.
\end{enumerate}
\mn
[Why?  As there is a model $M$ of cardinality $2^\sigma$ with $|\tau_M|
= \theta$ omitting a quantifier free type $p$ such that $M \subset N
\wedge M \equiv N \Rightarrow N$ realizes $p$.  Such $M$ exists as
$\sigma = \theta^{\aleph_0}$ and clause (b) of the assumption
$\circledast_1$ of \ref{e2}, \ref{e1}.]

Note
\mn
\begin{enumerate}
\item[$(*)_{13}$]   if $(a_{f,\alpha,u_1},a_{f,\alpha,u_2})$ is
$E_{2,\alpha}$-equivalent to $(a_{f,\alpha,v_1},a_{f,\alpha,v_2})$
then $G \models ``u_1-u_2 = v_1-v_2"$.
\end{enumerate}
\mn
[Why?  By induction on the $k$ from $(*)_4$.]

So to finish we assume toward contradiction
\mn
\begin{enumerate}
\item[$\boxtimes$]   $h$ is an isomorphism from $M_{1,\kappa}$
onto $M_{2,\kappa}$ which respects $(\pi_{1,\alpha},\pi_{2,\alpha})$
for $\alpha < \kappa$.
\end{enumerate}
\mn
So trivially
\mn
\begin{enumerate}
\item[$\circledast_1$]   $h(a_{f,\alpha,u}) \in \{a_{f,\alpha,v}:v \in G\}$
and $\bar a \in {}^n(A_\alpha) \Rightarrow h(\bar a) \in \bar
a/E_{n,\alpha} \subseteq \bar a/E'_{n,\alpha}$.
\end{enumerate}
\mn
[Why?  As $h$ respect $(\pi_{1,\kappa},\pi_{2,\kappa})$ see $(*)_1(e)$
and $(*)_{10}$ clearly $h(\bar a) \in \bar a/E'_{n,\alpha}$.  But
$h$ is an isomorphism from $M_{1,\kappa}$ onto $M_{2,\kappa}$ hence by
$(*)_6(b)$ we have $h(\bar a) \in (\bar a / E_{n,\alpha})$.]
\mn
\begin{enumerate}
\item[$\circledast_2$]   for $f \in {}^\kappa \sigma$ and $\alpha <
\kappa$ let $u_{f,\alpha} \in G$ be the $u \in G$ such that
$h(a_{f,\alpha,\emptyset}) = a_{f,\alpha,u}$
\sn
\item[$\circledast_3$]   for $f \in {}^\kappa \sigma,\alpha < \kappa$ and
$v \in G$ we have $h(a_{f,\alpha,v}) = a_{f,\alpha,v +_G u_{f,\alpha}}$.
\end{enumerate}
\mn
[Why?  By $\circledast_1$ clearly $h$ maps any finite 
sequence $\bar b \in {}^n(A_{1,\kappa})$ to an
$E_{n,\alpha}$-equivalent sequence for each $\alpha < \kappa$.  Now 
apply this to the pair $(a_{f,\alpha,\emptyset},a_{f,\alpha,u})$
recalling $(*)_{13}$.]
\mn
\begin{enumerate}
\item[$\circledast_4$]   we define a partial order $\le$ on ${}^\kappa
 \sigma$ as follows:

$f_1 \le f_2$ \Iff \,  there is a function $e \in {}^\sigma \sigma$ 
witnessing it; which means $f_1 = e \circ f_2$
\sn
\item[$\circledast_5$]  if $\alpha_1,\alpha_2 < \kappa$ and
$f_1 \le f_2$ (are from ${}^\kappa \sigma$) then $|u_{f_1,\alpha_1}| 
\le |u_{f_2,\alpha_2}|$.
\end{enumerate}
\mn
[Why?  This follows from $\circledast_6$ below.]
\mn
\begin{enumerate}
\item[$\circledast_6$]   if $e \in {}^\sigma \sigma,f_2 \in {}^\kappa
\theta$ and $f_1 = e \circ f_2 \in {}^\kappa \sigma$ and
$\alpha_1,\alpha_2 < \kappa$ then 

$u_{f_1,\alpha_1} \subseteq \{e(i):i \in u_{f_2,\alpha_2}\}$.
\end{enumerate}
\mn
[Why?  Choose $\alpha < \kappa$ such that $\alpha > \alpha_1,\alpha >
\alpha_2$ so $a_{f_1,\alpha_1,\emptyset},a_{f_2,\alpha_1,\emptyset}
\in M_{\ell,\alpha}$ for $\ell =1,2$.  
Recall that $h$ maps $R^{M_{1,\alpha}}_e$ onto
$R^{M_{2,\alpha}}_e$ by $\boxtimes$ and $R^{M_{2,\alpha}}_e =
R^{M_{1,\alpha}}_e$ because $g_\alpha$ maps $R^{M_{1,\alpha}}_e$ onto
itself (see the proof of $(*)_6$ above, the ``first" in that proof).   
Now see $(*)_6(d)$, i.e. the
definition of $R^{M_{1,\alpha}}_e$, i.e. obviously
$(a_{f_1,\alpha_1,\emptyset},a_{f_2,\alpha_2,\emptyset}) \in
R^{M_{1,\alpha}}_e$ so as $h$ is an isomorphism we have
$(h(a_{f_1,\alpha_1,\emptyset}),h(a_{f_2,\alpha_2,\emptyset})) \in
R^{M_{2,\alpha}}_e$ so by the previous sentence and the definitions of
$u_{f_\ell,\alpha_\ell}(\ell =1,2)$ in $\circledast_2$ we have
$(a_{f_1,\alpha_1,u_{f_1,\alpha_1}},a_{f_2,\alpha_2,u_{f_2,\alpha_2}}) \in
R^{M_{1,\alpha}}_e$ which by the definitions of $R^{M_{1,\alpha}}_e$
in $(*)_6(d)$ implies $u_{f_1,\alpha_1} \subseteq \{e(i):i \in
u_{f_2,\alpha_2}\}$ as promised.]
\mn
\begin{enumerate}
\item[$\circledast_7$]  $(a) \quad |u_{f,\alpha_1}| =
|u_{f,\alpha_2}|$ for $\alpha_1,\alpha_2 < \kappa,f \in {}^\kappa \sigma$
\sn
\item[${{}}$]    $(b) \quad \bold n(f) = |u_{f,\alpha}|$ is well  defined
\sn
\item[${{}}$]   $(c) \quad$ if $f_1 \le f_2$ then 
$\bold n(f_1) \le \bold n(f_2)$.
\end{enumerate}
\mn
[Why?  For clause (a) use $\circledast_6$ for the function 
$e = \text{ id}_\sigma$
and $f_1 = f_2 = f$.  Clause (b) follows.  Clause (c) holds by 
$\circledast_6$ equivalently by $\circledast_5$.]
\mn
\begin{enumerate}
\item[$\circledast_8$]  there are $f_* \in {}^\kappa \sigma$ and $\alpha_* <
\kappa$ such that:
\begin{enumerate}
\item[$(i)$]   if $f_* \le f \in {}^\kappa \sigma$ and $\alpha <
\kappa$ then $|u_{f_*,\alpha_*}| = |u_{f,\alpha}|$
\sn
\item[$(ii)$]   moreover if $f_* = e \circ f$ where $e \in {}^\sigma
\sigma$ and $f \in {}^\kappa \sigma,\alpha < \kappa$
\then \,  $e \restriction u_{f,\alpha}$ is one-to-one from $u_{f,\alpha}$ onto
$u_{f_*,\alpha}$ so $\bold n(f_*) = \bold n(f)$
\sn
\item[$(iii)$]   if $\alpha < \kappa,f_1 = e \circ f_2,
f_* = e_1 \circ f_1,f_* = e_2 \circ f_2$ so $e,e_1,e_2 \in {}^\sigma
\sigma$, \then \, $e \rest u_{f_2,\alpha}$ is one-to-one onto $u_{f_1,\alpha}$.
\end{enumerate}
\end{enumerate}
\mn
[Why?  First note that clause (ii), (iii) follows from clause (i).
Second, if clause (i) fails, then we can find a sequence $\langle
(f_n,\alpha_n,e_n):n < \omega \rangle$ such that
\mn
\begin{enumerate}
\item[$(\alpha)$]  $\alpha_n < \kappa,f_n \in {}^\kappa \sigma$ 
for $n < \omega$
\sn
\item[$(\beta)$]   $f_n \le f_{n+1}$ say $f_n = e_n
\circ f_{n+1}$ and $e_n \in {}^\sigma \sigma$
\sn
\item[$(\gamma)$]  $(e_n,f_{n+1},\alpha_{n+1})$ witness that
$(f_n,\alpha_n)$ does not satisfy the demand (i) on
$(f_*,\alpha_*)$ hence $\bold n(f_n) < \bold n(f_{n+1})$.
\end{enumerate}
\mn
Let $u_n = u_{f_n,\alpha_n}$ for $n < \omega$.  For $n < \omega$ and
$i < \sigma$ let $A_{n,i} = \langle \alpha <
\kappa:f_n(\alpha)=i\rangle$, so $\langle A_{n,i}:i < \sigma\rangle$
is a partition of $\kappa$ and $\alpha \in A_{n+1,i} \Rightarrow
\alpha \in A_{n,e_n(i)}$.  So letting $A_\eta = \cap\{A_{n,\eta(n)}:n
< \omega\}$ for $\eta \in {}^\omega \sigma$ clearly $\langle
A_\eta:\eta \in {}^\sigma \sigma\rangle$ is a partition of $\kappa$.

As we have $\sigma 
= \sigma^{\aleph_0}$ by $(*)_0$, there is a sequence $\langle
e^n:n < \omega\rangle$ satisfying 
$e^n \in {}^\sigma \sigma$ and $f \in {}^\kappa
\sigma$ such that $f_n = e^n \circ f$ for 
each $n < \omega$.  
So $n < \omega \Rightarrow f_n \le f$ which by $\circledast_7(c)$
implies $\bold n(f_n) \le \bold n(f)$.
As $\langle \bold n(f_n):n < \omega\rangle$ is increasing, 
easily we get a contradiction.]
\mn
\begin{enumerate}
\item[$\circledast_9$]   $\bold n(f_*) > 0$, i.e. $\alpha < \kappa
\Rightarrow u_{f_*,\alpha} \ne \emptyset$.
\end{enumerate}
\mn
[Why?  If $(\forall f \in {}^\kappa \sigma)(\forall \alpha < \kappa)
(u_{f,\alpha} = \emptyset)$ \then \, (by $\circledast_3$) we deduce 
$h$ is the identity, contradiction.
Otherwise assume $u_{f,\alpha} \ne \emptyset$ hence as in the proof of
$\circledast_8$ there is $f'$ such that $f_* \le f' \wedge f \le f'$ so by 
$\circledast_5$ and $\circledast_8$ we have 
$0 < |u_{f,\alpha}| \le |u_{f',\alpha}| = |u_{f_*,\alpha_*}|$.]
\mn
\begin{enumerate}
\item[$\circledast_{10}$]    if $f \in {}^\kappa \sigma,\alpha < \kappa$ and
$i \in u_{f,\alpha}$ then $\kappa = \sup\{\beta < \kappa:\alpha <
\beta$ and $f(\beta) =i\}$.
\end{enumerate}
\mn
[Why?  If not, let $\beta(*) < \kappa$ be $> \sup\{\beta <
\kappa:\alpha < \beta,f(\beta)=i\}$ and $> \omega$ and $> \alpha$.

\noindent
Let $Y = \{(a_{f,\alpha,u},a_{f,\beta(*),u}):u \in G,i \notin u\}$.
Now for every $\beta \in (\beta(*),\kappa)$ the function 
$g_\beta$ maps the set $Y$ onto itself (see its definition in $(*)_2$(f))
hence by the definition of $E_{2,\beta(*)+1}$ (in $(*)_4$) 
it follows that $\bar a
\in Y \Rightarrow \bar a/E_{2,\beta(*)+1} \subseteq Y$ and as $h$
respects $(\pi_{1,\beta(*)+1},\pi_{2,\beta(*)+1})$ it follows that
$h(\bar a) \subseteq \bar a/E'_{2,\beta(*)+1}$ and so $\kappa > \gamma
\ge \beta(*)+1 \Rightarrow g^{-1}_\gamma(h(\bar a)) \in \bar
a/E'_{2,\beta(*)+1}$.

Now for $\bar a \in Y$, the pairs $\bar a,h(\bar a)$ realizes the same
quantifier free type in $M_{1,\beta(*)+1},M_{2,\beta(*)+1}$
respectively, hence by the choice of $M_{2,\beta(*)+1}$ the pairs
$\bar a,g^{-1}_{\beta(*)+2} h(\bar a))$ realize the same quantifier
free type in $M_{1,\alpha}$.  By $(*)_6(b)$ recalling
$g^{-1}_{\beta(*)+2}(h(\bar a)) \in \bar a/E'_{2,\beta(*)+1}$ this
implies that $\bar a,g^{-1}_{\beta(*)+2}(h(\bar a))$ are
$E_{2,\beta(*)+1}$-equivalent.  By the definition of
$E_{2,\beta(*)+1},g^{-1}_{\beta(*)+2}(h(\bar a))$ belongs to the closure
of $\{\bar a\}$ under $\{g^{\pm 1}_\gamma:\gamma \in
(\beta(*);\kappa)\}$ hence $h(\bar a)$ belongs to it.  But by an earlier
sentence $Y$ is closed under those functions so $h(\bar a) \in Y$.
Similarly $h^{-1}(\bar a) \in Y$, hence $h$ maps $Y$ onto itself, 
recalling $\circledast_2$ this implies $i \notin u_{f,\alpha}$,
contradicting an assumption of $\circledast_{10}$, so
$\circledast_{10}$ holds.]
\smallskip

Now fix $f_*,\alpha_*$ as in $\circledast_8$ 
for the rest of the proof, \wilog \, $f_*$ is
onto $\sigma$ and let $u_{f_*,\alpha_*} = \{i^*_\ell:i <
\ell(*)\}$ with $\langle i^*_\ell:\ell < \ell(*)\rangle$ increasing for
simplicity.  Now for every $f \in {}^\kappa \sigma$ such that $f_* \le
f$ and $\alpha < \kappa$  by $\circledast_8(ii),(iii)$ 
we know that if $e \in {}^\sigma \sigma \wedge f_* = e \circ f$ then $e$ is
a one-to-one mapping from $u_{f,\alpha}$ onto $u_{f_*,\alpha_*}$; but so
$e \restriction u_{f,\alpha}$ is uniquely determined by
$(f_*,\alpha_*,f,\alpha)$ so let
$i_{f,\alpha,\ell} \in u_{f,\alpha}$ be the unique $i \in
u_{f,\alpha}$ such that $e(i) = i^*_\ell$ (equivalently $(\exists
\alpha)(f(\alpha) = i \wedge f_*(\alpha) = i^*_\ell$)).

Now if $f_* \le f \in {}^\kappa \sigma$ and 
$\alpha_1,\alpha_2 < \kappa$ and 
we choose $e = \text{ id}_\sigma$ so necessarily
$f \restriction u_{f,\alpha_1} = 
e \circ f \restriction u_{f,\alpha_2}$, \then \, 
$e \restriction \text{\rm Rang}(f \restriction
u_{f,\alpha_2})$ map $u_{f,\alpha_2}$ onto $u_{f,\alpha_1}$ but $e$
is the identity so we can write $u_f$ instead of
$u_{f,\alpha}$ let $i_{f,\ell} = i_{f,\alpha,\ell}$ for $\ell <
\ell(*),\alpha < \kappa$.

Let

\[
{\cA} = \{A \subseteq \kappa:\text{ for some } f,f_* \le f 
\text{ and } \alpha < \kappa \text{ we have } 
f^{-1}\{i_{f,0}\} \backslash \alpha \subseteq A\}
\]
\mn
\begin{enumerate}
\item[$\boxdot_1$]   ${\cA} \subseteq {\cP}(\kappa) \backslash
[\kappa]^{< \kappa}$.
\end{enumerate}
\mn
[Why?  As $\kappa$ is regular, this means $A \in {\cA} \Rightarrow 
A \subseteq \kappa \wedge \sup(A) = \kappa$ which holds by $\circledast_{10}$.]
\mn
\begin{enumerate}
\item[$\boxdot_2$]  $\kappa \in {\cA}$.
\end{enumerate}
\mn
[Why?  By the definition of ${\cA}$.]
\mn
\begin{enumerate}
\item[$\boxdot_3$]   if $A \in {\cA}$ and $A \subseteq B \subseteq
\kappa$ then $B \in \kappa$.
\end{enumerate}
\mn
[Why?  By the definition of ${\cA}$.]
\mn
\begin{enumerate}
\item[$\boxdot_4$]   if $A_1,A_2 \in {\cA}$ then $A =: A_1 \cap
A_2$ belongs to ${\cA}$.
\end{enumerate}
\mn
[Why?  Let $(f_\ell,e_\ell,\alpha_\ell)$ be such that
$f_* = e_\ell \circ f_\ell$ and $f_\ell \in {}^\kappa
\sigma,\alpha_\ell < \kappa$ and $f^{-1}_\ell \{i_{f_\ell,0}\}
\backslash \alpha_\ell \subseteq A_\ell$ for $\ell=1,2$.  Let pr:$\sigma \times
\sigma \rightarrow \sigma$ be one-to-one and onto and define $f \in {}^\kappa
\sigma$ by $f(\alpha) = \text{\rm pr}(f_1(\alpha),f_2(\alpha))$.
Clearly $f_\ell \le f$ for $\ell=1,2$ 
hence $i_{f,0}$ is well defined and $i_{f,0} =
\text{\rm pr}(i_{f_1,0},i_{f_2,0})$.  Now for every $\alpha < \kappa,
f(\alpha) = i_{f,0} \Rightarrow f_1(\alpha) = i_{f_1,0} 
\wedge f_2(\alpha) = i \Rightarrow
\alpha \in A_1 \wedge \alpha \in A_2 \Rightarrow \alpha \in A_1 \cap
A_2 \Rightarrow \alpha \in A$ so $f^{-1}\{i_{f,0}\} \subseteq A$ hence
$A \in \cA$.]
\mn
\begin{enumerate}
\item[$\boxdot_5$]   if $A \subseteq \kappa$ then $A \in {\cA}$
or $\kappa \backslash A \in {\cA}$.
\end{enumerate}
\mn
[Why?  Define $f \in {}^\kappa \sigma$:

\[
f(\alpha) = \begin{cases} 2f_*(\alpha) \quad &\text{ if } \alpha \in A \\
  2f_*(\alpha)+1 \quad &\text{ if } \alpha \in \kappa \backslash A.
\end{cases}
\]

\mn
Let $i = i_{f,0}$ so by the definition of $\cA$ we have $f^{-1}\{i\} =
f^{-1}\{i_{f,0}\} \in \cA$. But if $i$ is even then 
$f^{-1}\{i\} \subseteq A$ and $i$ is odd
then $f^{-1}\{i\} \subseteq \kappa \backslash A$ so by $\boxdot_3$ we
are done.]
\mn
\begin{enumerate}
\item[$\boxdot_6$]   ${\cA}$ is a uniform ultrafilter on $\kappa$.
\end{enumerate}
\mn
[Why?  By $\boxdot_1 - \boxdot_5$.]
\mn
\begin{enumerate}
\item[$\boxdot_7$]    ${\cA}$ is $\sigma^+$-complete.
\end{enumerate}
\mn
[Why?  Assume $B_\varepsilon \in {\cA}$ for $\varepsilon < \sigma$ and
  let $B = \cap\{B_\varepsilon:\varepsilon < \sigma\}$.
Define $A_\varepsilon \subseteq \kappa$ for $\varepsilon < \sigma$ as
follows: $A_{1 + \varepsilon} = \bigcap\limits_{\zeta < \varepsilon}
B_\zeta \backslash B_\zeta$ (so is $\kappa \backslash B_0$ if
$\varepsilon =0$) for $\varepsilon < \sigma$ and $A_0 = B$.
 Clearly $\langle A_\varepsilon:
\varepsilon <  \sigma \rangle$ is a partition of
$\kappa$, let $f \in {}^\kappa \sigma$ be such that $f \restriction
A_\varepsilon$ is constantly $\varepsilon$.  Let $f' \in {}^\kappa
\theta$ be such that $f \le f' \wedge f_* \le f'$.  Now
$(f')^{-1}\{i_{f',0}\} \in {\cA}$ is included in some
$A_\varepsilon$. If $\varepsilon =0$ this exemplifies
$\bigcap\limits_{\varepsilon < \sigma} B_\varepsilon \in {\cA}$ as required.
If $\varepsilon = 1 + \zeta < \sigma$, then $(f')^{-1}\{i_{f',0}\}
\subseteq A_\varepsilon \subseteq \kappa \backslash B_\varepsilon$,
contradiction to $\boxdot_6$ because $B_\varepsilon \in {\cA}$ and
$(f')^{-1}\{i_{f',0}\} \in {\cA}$.]

So by the assumptions of \ref{e2}, that is, $\circledast_1(b)$ of
\ref{e1} we get a contradiction, coming from the assumption ``toward
contradiction (j) of $\boxplus$ of \ref{e2} fails", so it holds and
the other clauses were proved so we are done.  
\end{PROOF}

\begin{theorem}
\label{e3}  For every $\theta$ there is an ${\gk} = {\gk}^*_\theta$ such that
\mn
\begin{enumerate}
\item[$\otimes$]  $(a) \quad {\gk}$ is an a.e.c. with 
{\rm LST}$({\gk}) = \theta,|\tau_{\gk}| = \theta$
\sn
\item[${{}}$]   $(b) \quad {\gk}$ has the amalgamation property
\sn
\item[${{}}$]  $(c) \quad {\gk}$ admits intersections (see Definition
  \ref{e4} below) 
\sn
\item[${{}}$]   $(d) \quad$ if $\kappa$ is a regular cardinal and
there is no uniform $\theta^+$-complete 

\hskip25pt ultrafilter  on $\kappa$, \then \,: ${\gk}$ is not 
$(\le 2^\kappa,\kappa)$-sequence-local for types,  

\hskip25pt i.e., we can find an $\le_{\gk}$-increasing
continuous sequence $\langle M_i:i \le \kappa \rangle$  

\hskip25pt of models and $p \ne q \in \mathscr{S}_{\gk}(M_\kappa)$
such that $i < \kappa \Rightarrow p \restriction M_i = q \restriction M_i$

\hskip25pt  and $M_\kappa$ is of cardinality $\le 2^\kappa$.
\end{enumerate}
\end{theorem}

We shall prove \ref{e3} below.  As in \cite[1.2,\S4]{BlSh:862} the aim of the 
definition of ``admit intersections" is to ensure types behave reasonably.
\begin{definition}
\label{e4}  
We say an a.e.c. ${\gk}$ admits intersections \when \, there is a
function $c \ell_{\gk}$ such that:
\mn
\begin{enumerate}
\item[$(a)$]  $c \ell_{\gk}(A,M)$ is well defined iff $M \in K_{\gk}$
  and $A \subseteq M$
\sn
\item[$(b)$]  $c \ell_{\gk}(A,M)$ is preserved under isomorphisms and 
$\le_{\gk}$-extensions
\sn
\item[$(c)$]  for every $M \in K_{\gk}$ and non-empty $A
\subseteq M$ the set is $B = c \ell_{\gk}(A,M)$ satisfies: $M \rest
B \in K_{\gk},M \rest B \le_{\gk} M$ and 
$A \subseteq M_1 \le_{\gk} N \wedge M \le_{\gk} N \Rightarrow B
\subseteq M_1$; we may  use $cl_{\gk}(A,M)$ for $M \rest c
\ell_{\gk}(A,M)$.
\end{enumerate}
\end{definition}

\begin{claim}
\label{e5}  
Assume ${\gk}$ is an {\rm a.e.c.} admitting intersections.  \Then \,
{\bf{\rm tp}}$_{\gk}(a_1,M,N_1) = \text{\rm tp}_{\gk}(a_2,M,N_2)$ iff 
letting $M_\ell = N_\ell \rest c
\ell_{\gk}(M \cup \{a_\ell\})$, there is an isomorphism from $M_1$
onto $M_2$ over $M$ mapping $a_1$ to $a_2$.
\end{claim}

\begin{PROOF}{\ref{e5}}
Should be clear by the definition.
\end{PROOF}

\begin{remark}
\label{e5d}
In Theorem \ref{e3} we can many 
times demand $\|M_\kappa\| = \kappa$, e.g., if $(\exists \lambda)
(\kappa = 2^\lambda)$.
\end{remark}

Note we now show that \ref{e3} is best possible.
\begin{claim}
\label{e6}  
1) If ${\gk}$ satisfies clause (a) of \ref{e3}, (i.e. $\gk$ is an a.e.c. 
with $\LST$-number $\le \theta$ and $|\tau_{\gk}| \le \theta$)
and $\kappa$ fails the assumption of clause (d) of \ref{e3}, that
is there is a uniform $\theta^+$-complete ultrafilter on $\kappa$, 
\then \, the conclusion of clause (d) of \ref{e3} fails, that is
${\gk}$ is $\kappa$-sequence local for types.

\noindent
2) If $D$ is a $\theta^+$-complete ultrafilter on $\kappa$ and 
${\gk}$ is an a.e.c. with {\rm LST}$({\gk}) \le \theta$ \then \,
ultraproducts by $D$ preserve $``M \in {\gk}",``M \le_{\gk} N"$, i.e.
\mn
\begin{enumerate}
\item[$\boxtimes$]   if $M_i,N_i (i < \kappa)$ are
$\tau({\gK})$-models and $M = \prod\limits_{i < \kappa} M_i/D$ and $N = 
\prod\limits_{i < \kappa} N_i$ \then \,:
\begin{enumerate}
\item[$(a)$]   $M \in K$ if $\{i < \kappa:M_i \in {\gk}\} \in D$
\sn
\item[$(b)$]   $M \le_{\gk} N$ if $\{i:M_i \le_{\gk} N_i\} \in D$.
\end{enumerate}
\end{enumerate}
\end{claim}

\begin{PROOF}{\ref{e6}}
Note that if $D$ is $\theta^+$-complete, \then \, it is
$\sigma^+$-complete where $\sigma = \theta^{\aleph_0}$
(and much more, it is $\theta'$-complete for the first 
measurable $\theta' > \theta$).

\noindent
1) So assume
\mn
\begin{enumerate}
\item[$\boxplus$] 
\begin{enumerate}
\item[(a)]  $\langle M_i:i \le \kappa\rangle$ is
$\le_{\gk}$-increasing
\sn
\item[(b)]  $M_\kappa = N_0 \le_{\gk} N_\ell$ for $\ell=1,2$
\sn
\item[(c)]  $p_\ell =$ {\bf tp}$_{\gk}(a_\ell,N_0,N_\ell)$ for $\ell=1,2$
\sn
\item[(d)]  $i < \kappa \Rightarrow p_1 \rest M_i = p_2 \rest M_i$.
\end{enumerate}
\end{enumerate}
\mn
We shall show $p_1 = p_2$, this is enough.

Without loss of generality 
\mn
\begin{enumerate}
\item[$(*)_1$]  
\begin{enumerate}
\item[(a)]  $a_1 = a_2$ call it a 
\sn
\item[(b)]  $\tau_{\gk} \subseteq \cH(\theta)$.
\end{enumerate}
\end{enumerate}
\mn
By $(d)$ of $\boxplus$ we have:
\mn
\begin{enumerate}
\item[$(d)^+$]  for each $i < \kappa$ there are $n_i < \omega$ and $\langle
N_{i,m}:n \le n_i \rangle$ such that
\sn
\begin{enumerate}
\item[$(\alpha)$]  $N_{i,0} = N_1$
\sn
\item[$(\beta)$]  $N_{i,m_i} = N_2$ or just $h_i$ is an isomorphism
  from $N_{i,m_i}$ onto $N_2$ such that $h_i \rest (M_i \cup \{a\})$
  is the identity
\sn
\item[$(\gamma)$]  $a \in N_{i,\ell}$ and $M_i \le_{\gk} N_{i,\ell}$
\sn
\item[$(\delta)$]  if $m < m_i$ then $N_{i,2m +1} \le_{\gk}
N_{i,2m},N_{i,2m+2}$.
\end{enumerate}
\end{enumerate}
\mn
As $\kappa = \text{ cf}(\kappa) > \aleph_0$ \wilog \, $i < \kappa
\Rightarrow n_i = n_*$.
\mn
Let $\chi$ be such that $\langle M_i:i \le \kappa\rangle,\left<
\langle N_{i,n}:n \le n_*\rangle:i < \kappa\right>$ and
$\gk_{\text{LST}(\gk)}$ all belongs to $\cH(\chi)$; concerning
$\gk_{\LST(\gk)}$ this means $\tau_{\chi}$ and $\LST(\gk)$ belongs to
$\cH(\chi)$ hence $\{M \in K_{\gk}:M \in \cH(\LST^+_{\gk})\}$ and
$\le_{\gk} \rest \cH(\LST^+_{\gk})$ belongs to 
$\cH(\chi)$; those hold by $(*)_1(b)$.  Let $\gB$ be the
ultrapower $(\cH(\chi),\in)^\kappa/D$ and $\bold j_0$ the canonical
embedding of $(\cH(\chi),\in)$ into $\gB$ and let $\bold j_1$ be the
Mostowski-Collapse of $\gB$ to a transitive set $\cH$ and let $\bold
j = \bold j_1 \circ \bold j_0$.  So $\bold j$ is an elementary
embedding of $(\cH(\chi),\in)$ into $(\cH,\in)$ and even an
$\bbL_{\theta^+,\theta^+}$-elementary one.  Recall we are assuming
\wilog \, $\tau_{\gk} \subseteq \cH(\theta)$ hence $\bold j(\tau_{\gk}) =
\tau_{\gk}$ hence by part (2), $\bold j$ preserves $``N \in
K_{\gk}",``N^1 \le_{\gk} N^2"$.  ``$h$ is an isomorphism from $N'$
onto $N''$.  

So $\bold j(\langle M_i:i \le \kappa\rangle)$ has the form $\langle
M^*_i:i \le \bold j(\kappa)\rangle$ but $\bold j(\kappa) > \kappa_* :=
\bigcup\limits_{i <\kappa} \bold j(i)$ by the uniformity of $D$ and let
$\bold j(\left< \langle N_{i,n}:n \le n_*\rangle:i < \kappa \right> =
\left<\langle N^*_{i,n}:n \le n^*\rangle:i < \bold j(\kappa)\right>$
and $\bold j(\langle h_i:i < \kappa \rangle) = 
\langle h^*_i:i < \kappa_*\rangle$.

So
\mn
\begin{enumerate}
\item[(a)]  $\bold j \rest M_\kappa$ is a $\le_{\gk}$-embedding
of $M_\kappa$ into $M^*_\kappa$ hence even into $M^*_{\kappa_*}$
\sn
\item[(b)]  $M^*_{\kappa_*} \le_{\gk} N^*_{i,n}$ and $\bold j(a) \in
N^*_{i,n}$ for $i < \kappa,n \le n_*$
\sn
\item[(c)]  $N^*_{i,0} = \bold j(N_1)$
\sn
\item[(d)]  $h_{\kappa_*}$ is an isomorphism from $N_{\kappa_*,n_*}$
  onto $\bold j(N_2)$
\sn
\item[(e)]  $N^*_{\kappa_*,2m+1} \le_{\gk} N^*_{\kappa_*,2m},
N^*_{\kappa_*,2m+2}$ for $2m +1 < n_*$
\sn
\item[(f)]  $\bold j(a) \in N_{\kappa_*,m}$.
\end{enumerate}
\mn
Together we are done.

\noindent
2) By the representation theorem of a.e.c. \cite[\S1]{Sh:88r}.
\end{PROOF}

\begin{PROOF}{\ref{e3}}
\underline{Proof of \ref{e3}}

Let $\sigma = \theta^{\aleph_0}$.  Let $G =
([\sigma]^{<\aleph_0},\Delta)$ and let $\langle c_i:i < \sigma\rangle$
list the members of $G$, let $\langle \eta_\alpha:\alpha <
\sigma\rangle$ list ${}^\omega \theta$.

Now
\mn
\begin{enumerate}
\item[$\boxtimes_1$]  let $B_{\varepsilon,n} \subseteq G$ for $\varepsilon
< \theta$ be such that: if $a,b \in G$ then
$(\forall \varepsilon < \theta)(\forall n < \omega)
(a \in B_{\varepsilon,n} \equiv b \in B_{\varepsilon,n}) \Rightarrow
a=b$; moreover, $B_{\varepsilon,n} = \{c_\alpha:\eta_\alpha(n) =
\varepsilon\}$. 
\end{enumerate}
\mn
Let $\tau$ have the predicates $G,I,J,H$(unary), $E_1,Q$(binary), 
$R_{n,\alpha}$ ($n$-place; $n < \omega,\alpha < \theta),
P_{\varepsilon,n}$(unary; $\varepsilon < \theta)$ and 
function symbols $F_1$(unary), $F_2,\pi,+$(binary); 
so $|\tau| = \theta$.  We define $K$ as a class of $\tau$-models by:
\mn
\begin{enumerate}
\item[$\boxtimes_2$]   $M \in K$ \Iff \, (up to isomorphism):
\sn
\begin{enumerate}
\item[$(a)$]  $\langle G^M,I^M,J^M,H^M\rangle$ is a partition of
  $|M|$, (recall that they are unary)
\sn
\item[$(b)$]  $(G^M,+^M)$ is a subgroup of the group $([\sigma]^{< \aleph_0},
\Delta)$, $P^M_{\varepsilon,n} \subseteq G^M$ for $\varepsilon <
\theta,\langle P^M_{\varepsilon,n}:\varepsilon < \theta\rangle$ be a
partition of $M$ such that $a \ne b \in G^M \Rightarrow (\exists
\varepsilon < \theta)(\exists n < \omega)[a \in P^M_{\varepsilon,n}
\wedge b \notin P^M_{\varepsilon,n}]$
\sn
\item[$(c)$]  $Q^M \subseteq H^M \times J^M$ is such that $(\forall a
  \in H^M)(\exists^{\le 1} b)((a,b) \in Q^M)$;
\sn
\item[$(d)$]  $E^M_1$ is an equivalence relation on $H^M$ such that:
  if $(a_1,b) \in Q^M$ and $a_2 \in H^M$ then $a_1 E^M_1 a_2
  \Leftrightarrow (a_2,b) \in Q^M$
\sn
\item[$(e)$]  $\pi^M$ is a function from $H^M$ into $I^M$
\sn
\item[$(f)$]  $E^M_2 = \{(a,b):a E^M_1 b$ and $\pi^M(a) = \pi^M(b)$ so
  $a,b \in H\}$
\sn
\item[$(g)$]  $F^M_2$ is a partial two-place function such that:
\sn
\item[${{}}$]  $(\alpha) \quad F^M_2(a,b)$ is well defined iff
$b \in G^M,a \in H^M$
\sn
\item[${{}}$]  $(\beta) \quad$ for $a \in H^M,\langle
F^M_2(a,b):b \in G^M\rangle$ list $a/E^M_2 = \{a' \in H^M$:

\hskip35pt $\pi^M(a') = \pi^M(a)\}$ with no repetitions
\sn
\item[${{}}$]  $(\gamma) \quad$ if $a \in H^M$ and $b,c \in G$
then $F^M_2(a,b +^M c) = F^M_2(F^M_2(a,b),c)$,

\hskip35pt  on the + see clause (b)
\sn
\item[${{}}$]   $(\delta) \quad F(a,0_{G^M}) = a$ for $a \in H^M$
\sn
\item[${{}}$]  $(\varepsilon) \quad$ for $n < \omega$ and $\gamma <
\theta$ the relation $R^M_{n,\gamma}$ is an $n$-place relation 

\hskip35pt $\subseteq \cup\{{}^n(a/E^M_1):a \in H^M\}$.
\end{enumerate}
\end{enumerate}
\mn
We define $\le_{\frak k}$ as being a submodel.  Easily
\mn
\begin{enumerate}
\item[$\boxtimes_3$]  ${\gk} = (K,\le_{\gk})$ is an a.e.c.
\end{enumerate}
\mn
For $A \subseteq M \in K$ let
\mn
\begin{enumerate}
\item[$(a)$]   $c \ell^1_M(A) = \text{ the subgroup of } 
(G^M,+^M) \text{ generated by } (A \cap G^M) \cup 
\{b \in G^M:\text{ for some } a_1 \ne a_2 \in A \cap H^M$ 
we have $a_1 E^M_2 a_2 \text{ and } F^M_2(a_1,b)=a_2\}$
\sn
\item[$(b)$]   $c \ell^2_M(A) = (A \cap I^M) \cup \{\pi^M(a):
a \in A \cap H^M\}$
\sn
\item[$(c)$]   $c \ell^3_M(A) = \{a \in H^M:\text{ for some } b \in c
\ell^0_M(A) \text{ and } a'_1 \in A \cap H^M$ we have $a = F^M_2(a'_1,b)\}$
\sn
\item[$(d)$]   $c \ell(A,M) = c \ell_M(A) = M \restriction (\cup
\{c \ell^\ell_M(A):\ell =1,2,3\})$.
\end{enumerate}
\mn
Now this function $c \ell(A,M)$ shows that ${\gk}$ admits intersections (see
Definition \ref{e4}) so
\mn
\begin{enumerate}
\item[$\boxtimes_4$]   ${\gk}$ admits closure and LST$(\gk) +
|\tau_{\gk}| = \theta$.
\end{enumerate}
\mn
Assume $\kappa$ is as in clause (d) of \ref{e3}, 
we use the $M_{\ell,\alpha}(\ell=1,2,\alpha \le \kappa)$
constructed in \ref{e2} (the relevant properties are stated in \ref{e2}).
They are not in the right vocabulary so let $M'_{\ell,\alpha}$ be the
following $\tau$-model:
\mn
\begin{enumerate}
\item[$\boxtimes_5$]   $(a) \quad$ \underline{elements} 
$\quad G^{M'_{\ell,\alpha}} = G$

\hskip35pt $I^{M'_{\ell,\alpha}} = I_\alpha$

\hskip35pt $J^{M'_{\ell,\alpha}} = \{t^*_\ell\},t^*_\ell$ just a
new element

\hskip35pt $H^{M'_{\ell,\alpha}} = |M_{\ell,\alpha}|$

(we assume disjointness)
\sn
\item[${{}}$]  $(b) \quad (G^{M'_{\ell,\alpha}},
+^{M'_{\ell,\alpha}})$ is $G =([\sigma]^{<\aleph_0},\Delta)$

\hskip35pt $P^{M'_{\ell,\alpha}}_\varepsilon \subseteq G^M$ as required in
$\boxtimes_1$ \underline{not} depending on $(\ell,\alpha)$
\sn
\item[${{}}$]  $(c) \quad F^{M'_{\ell,\alpha}}_1$ is 
constantly $t^*_\ell$ on $H^{M'_{\ell,\alpha}}$
\sn
\item[${{}}$]   $(d) \quad E^{M'_{\ell,\alpha}}_1 =
\{(a,b):F^{M'_{\ell,\alpha}}_1(a) = F^{M'_{\ell,\alpha}}_1(b)$ so $a,b
\in H^{M'_{\ell,\alpha}}$
\sn
\item[${{}}$]   $(e) \quad \pi^{M'_{\ell,\alpha}}$ 
is $\pi_{\ell,\alpha}$ (constructed in \ref{e2})
\sn
\item[${{}}$]  $(f) \quad E^{M'_{\ell,\alpha}}_2 = \{(a,b):a
E^{M'_{\ell,\alpha}}_1 b$ and $\pi^{M'_{\ell,\alpha}}(a) = 
\pi^{M'_{\ell,\alpha}}(b)$ so $a,b \in H^{M'_{\ell,\alpha}}\}$
\sn
\item[${{}}$]   $(g) \quad F^{M'_{\ell,\alpha}}_2(a,b) =
F^{M_{\ell,\alpha}}_b(b)$ for $a \in H^{M'_{\ell,\alpha}}$
\sn
\item[${{}}$]   $(h) \quad R^{M'_{\ell,\alpha}}_{\gamma,n}$ for $n <
\omega,\gamma < \sigma$ list the relations of $M_{\ell,\alpha}$.
\end{enumerate}
\mn
Let $M'_{0,\alpha} = M'_{\ell,\alpha} \restriction 
(G^{M'_{\ell,\alpha}} \cup I^{M'_{\ell,\alpha}})$ for $\ell = 1,2$ and
$\alpha \le \kappa$ (we get the same result).

Note easily
\mn
\begin{enumerate}
\item[$\boxtimes_6$]   $M_{0,\alpha} \le_{\frak k}
M_{\ell,\alpha},\langle M_{\ell,\alpha}:\alpha \le \kappa \rangle$ is
$\le_{\gk}$-increasing (check)
\sn
\item[$\boxtimes_7$]  tp$_{\gk}(t^*_1,M'_{0,\alpha},M'_{1,\alpha})
= \text{\rm tp}_{\gk}(t^*_2,M'_{0,\alpha},M'_{2,\alpha})$ for
$\alpha < \kappa$.
\end{enumerate}
\mn
[Why?  By the isomorphism from $M_{1,\alpha}$ onto $M_{2,\alpha}$
respecting $(\pi_{1,\alpha},\pi_{2,\alpha})$ in \ref{e1}.]
\mn
\begin{enumerate}
\item[$\boxtimes_8$]   tp$_{\gk}(t^*_1,M'_{0,\kappa},M'_{1,\kappa})
\ne \text{\rm tp}_{\gk}(t^*_2,M'_{0,\kappa},M'_{2,\kappa})$.
\end{enumerate}
\mn
[Why?  By the non-isomorphism in \ref{e1}; extension will not help.]

Now by the``translation theorem" of \cite[4.7]{BlSh:862} 
we can find ${\gk}'$ which has all the needed properties,
i.e. also the amalgamation and JEP. 
\end{PROOF}
\newpage

\section {Compactness of types in a.e.c.} 

Baldwin \cite{Bal0x} ask ``can we in ZFC prove that some a.e.c. has
amalgamation, JEP but fail compactness of types".  The background is
that in \cite{BlSh:862}
we construct one using diamonds.

To me the question is to show this class can be very large (in ZFC).

Here we omit amalgamation and 
accomplish both by direct translations of problems of
existence of models for theories in $\bbL_{\kappa^+,\kappa^+}$, first
in the propositional logic.  So whereas in \cite{BlSh:862} we have an
original group $G^M$,
 here instead we have a set $P^M$ of propositional ``variables"
and $P^M$, set of such sentences (and relations and functions explicating
this; so really we use coding but are a little sloppy in stating this
obvious translation).

In \cite{BlSh:862} we have $I^M$, set of indexes, $0$ and $H$, set of
Whitehead cases, $H_t$ for $t \in I^M$, here we have $I^M$, each $t \in
I^N$ representing a theory $P^M_t \subseteq P^M$ and in $J^M$ we give
each $t \in I^M$ some models ${\cM}^M_s:P^M \rightarrow
\{\text{true,false}\}$.  This is set up so that amalgamation holds.

\begin{notation}  In this section types are denoted by $\bold p,\bold q$
as $p,q$ are used for propositional variables.
\end{notation}

\begin{definition}
\label{b2.0}  
1) We say that an a.e.c. ${\gk}$ has
$(\le \lambda,\kappa)$-sequence-compactness (for types) 
\when \,: if $\langle M_i:i
\le \kappa\rangle$ is $\le_{\gk}$-increasing continuous and $i <
\kappa \Rightarrow \|M_i\| \le \lambda$ and $\bold p_i \in 
\mathscr{S}^{< \omega}(M_i)$ for $i < \kappa$ satisfying 
$i < j < \kappa \Rightarrow \bold p_i = \bold p_j \rest M_i$ 
\then \, there is $\bold p_\kappa \in \mathscr{S}^{< \omega}(M_\kappa)$
such that $i < \kappa \Rightarrow p_\kappa \rest M_i = p_i$.

\noindent
2) We define ``$(=\lambda,\kappa)$-sequence-compactness" similarly.  Let
$(\lambda,\kappa)$-sequence-compactness mean 
$(\le \lambda,\kappa)$-compactness.
\end{definition}

\begin{question}
\label{2b.1}  Can we find an a.e.c. ${\gk}$ with
amalgamation and JEP such that $\{\theta:{\gk}$ have
$(\lambda,\theta)$-compactness of types for every $\lambda\}$ is
complicated, say:
\mn
\begin{enumerate}
\item[$(a)$]   not an end segment but with ``large" members
\sn
\item[$(b)$]   any $\{\theta:\theta$ satisfies $\psi\},\psi \in 
\bbL_{\kappa^+,\kappa^+}$ (second order).
\end{enumerate}
\end{question}

\begin{definition}
\label{b2.2}  Let $\kappa \ge \aleph_0$, we
define ${\gk} = {\gk}_\kappa$ as follows:
\mn
\begin{enumerate}
\item[$(A)$]   the vocabulary $\tau_{\frak k}$ consist of $F_i(i \le
\kappa),R_\ell(\ell=1,2),P,\Gamma,I,J,c_i \, (i < \kappa),F_i(i \le
\kappa)$, (pedantically see later),
\sn
\item[$(B)$]    the universe of $M \in K_{\gk}$ is the 
disjoint union of $P^M,\Gamma^M,I^M,J^M$ so $P,\Gamma,I,J$
are unary predicates
\sn
\item[$(C)$]   $(a) \quad P^M$ a set of propositional variables
(i.e. this is how we treat them)
\sn
\item[${{}}$]   $(b) \quad \Gamma^M$ is a set of sentences of one of
the forms $\varphi = (p),\varphi = (r \equiv p \wedge$

\hskip25pt $q),\varphi = 
(q \equiv \neg p),\varphi = (q \equiv \bigwedge\limits_{i <
\kappa} p_i)$, so $p,q,p_i \in P^M$ 

\hskip25pt but in the last case
$\{p_i:i < \kappa\} \subseteq \{c^M_i:i < \kappa\}$ (or code this!)
\sn
\item[${{}}$]   $(c) \quad$ for $i < \kappa$ the function 
$F^M_i:\Gamma^M \rightarrow P^M$ are such that for every

\hskip25pt  $i < \kappa$ and $\varphi \in \Gamma^M$ we have:
\sn
\item[${{}}$]  $\qquad (\alpha) \quad$ if $\varphi = (p)$ and $i \le
\kappa$ then $F_{1+i}(\varphi) = p,F_0(\varphi) = c_0$
\sn
\item[${{}}$]  $\qquad (\beta) \quad$ if $\varphi = (r \equiv p \wedge
q)$ \then \, $F_i(\varphi)$ is $c_1$ if $i=0$, is $p$ if $i=1$, is $q$

\hskip45pt if $r=2$ is $r$ if $r \ge 3$
\sn
\item[${{}}$]  $\qquad (\gamma) \quad$ if $\varphi = (q \equiv \neg p)$
then $F_i(\varphi)$ is $c_2$ if $i=0$, $p$ if $i=1$, $q$ if $i \ge 2$
\sn
\item[${{}}$]  $\qquad (\delta) \quad$ if $\varphi = (q
\equiv \bigwedge \limits_{j < \kappa} p_j)$ then $F_i(\varphi)$ is
$c_3$ if $i=0$,

\hskip45pt $q$ if $i=1,p_{2+j}$ if $i =j+1$
\sn
\item[${{}}$]   $(d) \quad I$ a set of theories, i.e. $R^M_1 \subseteq \Gamma
\times I$ and for $t \in I$ let 

\hskip25pt $\Gamma^M_t = \{\psi \in \Gamma^M:\psi R^M_1 t\} \subseteq \Gamma^M$
\sn
\item[${{}}$]   $(e) \quad J$ is a set of models, i.e. $R^M_2
\subseteq (\Gamma \cup P) \times J$ and for $s \in J$ we have 

\hskip25pt  ${\cM}^M_s$ is the model, i.e. function giving truth values to 
$p \in P^M$, i.e. 
\sn
\begin{enumerate}
\item[${{}}$]  $(\alpha) \quad {\cM}^M_s(p)$ is true if $p_i R^M_2 s$; 
is false if $\neg p R^M_2 s$
\sn
\item[${{}}$]  $(\beta) \quad (\varphi,s) \in R^M_2$ iff
computing the truth value of $\varphi$ in $\cM^M_s$ 

\hskip25pt we get truth
\end{enumerate}
\sn
\item[${{}}$]   $(f) \quad F^M_\kappa:J^M \rightarrow I^M$ such 
that $s \in J^M \Rightarrow {\cM}^M_s$ is a model of $\Gamma_{F^M_\kappa(s)}$
\sn
\item[${{}}$]  $(g) \quad (\forall t \in I^M)(\exists s \in
J^M)(F^M_\kappa(s)=t)$
\sn
\item[$(D)$]   $M \le_{\gk} N$ \Iff \,  $M \subseteq N$ are
$\tau_{\gk}$-models from $K_{\gk}$.
\end{enumerate}
\end{definition}

\begin{claim}
\label{b2.3}  
${\gk}$ is an a.e.c., LST$({\gk}) = \kappa$.
\end{claim}

\begin{PROOF}{\ref{b2.3}}
Obvious.  
\end{PROOF}

\begin{claim}
\label{b2.5}  ${\gk}$ has the JEP.
\end{claim}

\begin{PROOF}{\ref{b2.5}}
Just like disjoint unions (also of the relations and functions) except
for the individual constants $c_i$ (for $i < \kappa$).
\end{PROOF}

\begin{claim}
\label{b2.6}  
Assume $M_0 \le_{\gk} M_\ell$ for $\ell=0,1$ and $|M_0| = 
P^{M_0} \cup \Gamma^{M_0} = P^{M_\ell} \cup
\Gamma^{M_\ell}$ for $\ell=1,2$ and $a_\ell \in I^{M_\ell}$ for
$\ell=1,2$.  \Then \, $\text{\bf tp}_{\gk}(a_1,M_0,M_1) 
=$ {\bf tp}$_{\gk}(a_2,M_0,M_2)$ \Iff \, 
$\Gamma^{M_1}_{a_1} = \Gamma^{M_2}_{a_2}$.
\end{claim}

\begin{PROOF}{\ref{b2.6}}
\underline{The if direction, $\Leftarrow$}

Let $h$ be a one to one mapping with domain $M_1$ such that $h \rest
M_0 =$ the identity, $h(a_1) = a_2$ and $h(M_1) \cap M_2 = M_0 \cup
\{a_2\}$.  Renaming \wilog \, $h$ is the identity.  Now define $M_3$ as
$M_1 \cup M_2$, as in \ref{b2.5}, now $a_1 = a_2$ does not cause
trouble because $P^{M_0} = P^{M_\ell},\Gamma^{M_0} = \Gamma^{M_\ell}$
for $\ell=1,2$.
\bigskip

\noindent
\underline{The only if direction, $\Rightarrow$}

Obvious.  
\end{PROOF}

\begin{claim}
\label{b2.11}  Assume $\lambda,\theta$ are such that:
\mn
\begin{enumerate}
\item[$(a)$]   $\theta$ is regular $\le \lambda$ and $\lambda \ge \kappa$
\sn
\item[$(b)$]   $\langle \Gamma_i:i \le \theta \rangle$ is
$\subseteq$-increasing continuous sequence of sets propositional sentences
in $\bbL_{\kappa^+,\omega}$ such that [$\Gamma_i$ has a model
$\Leftrightarrow i < \theta$]
\sn
\item[$(c)$]   $|\Gamma_\theta| \le \lambda$.
\end{enumerate}
\mn
\Then \, ${\gk}$ fail $(\lambda,\theta)$-sequence-compactness (for types).
\end{claim}

\begin{remark}
We may wonder but: for $\theta = \aleph_0$ compactness holds?  Yes,
but only assuming amalgamation.
\end{remark}

\begin{PROOF}{\ref{b2.11}}
 Without loss of generality $|\Gamma_0| = \lambda$.
Without loss of generality $\langle p^*_\varepsilon:\varepsilon 
< \kappa\rangle$ are pairwise distinct
propositions variables appearing in $\Gamma_0$ (but not necessarily
$\in \Gamma_0$) and 
each $\psi \in \Gamma_i$ is of the form $(p)$ or
$r \equiv p \wedge q$ or $r \equiv \neg p$ or 
$r \equiv \bigwedge\limits_{i < \kappa} p_i$ where $\{p_i:i < \kappa\}
\subseteq \{p^*_\varepsilon:\varepsilon < \kappa\}$.

Let $P_i$ be the set of propositional variables appearing in $\Gamma_i$
without loss of generality $|P_i| = \lambda$.

We choose a model $M_i$ for $i \le \theta$ such that:
\mn
\begin{enumerate}
\item[$\boxplus$]   $(a) \quad |M_i| = P_i \cup \Gamma_i$
\sn
\item[${{}}$]   $(b) \quad P^M = P_i$ and $\Gamma^{M_i} = \Gamma_i$
\sn
\item[${{}}$]   $(c) \quad$ the natural relations and functions.
\end{enumerate}
\mn
Let ${\cM}_i:P_i \rightarrow \{$true false$\}$ be a model of
$\Gamma_i$.

We define a model $N_i \in K_{\gk}$ for $i < \kappa$ (but not for $i =
\theta$!)
\mn
\begin{enumerate}
\item[$\boxtimes$]   $(a) \quad M_i \le_{\gk} N_i$
\sn
\item[${{}}$]   $(b) \quad P^{N_i} = P^{M_i}$
\sn
\item[${{}}$]   $(c) \quad \Gamma^{N_i} = \Gamma^{M_i}$
\sn
\item[${{}}$]   $(d) \quad I^M = \{t_i\}$
\sn
\item[${{}}$]   $(e) \quad J^M = \{s_i\}$
\sn
\item[${{}}$]   $(f) \quad F^{N_i}_\kappa(s_i) = t_i$
\sn
\item[${{}}$]   $(g) \quad R^{N_i}_1 = \Gamma_i \times \{t_i\}$
\sn
\item[${{}}$]   $(h) \quad R^{N_i}_2$ is chosen such that 
${\cM}^{N_i}_{s_i}$ is ${\cM}_i$
\sn
\item[${{}}$]   $(i) \quad F^{N_i}_i \, (i < \kappa)$ are defined naturally.
\end{enumerate}
\mn
Now
\mn
\begin{enumerate}
\item[$(*)_1$]  $\bold p_i = \text{\bf tp}_{\gk}(t_i,M_i,N_i) \in \bold
S^1(M_i)$.
\end{enumerate}
\mn
[Why?  Trivial.]
\mn
\begin{enumerate}
\item[$(*)_2$]   $i<j < \theta \rightarrow \bold p_i = \bold p_j \rest M_j$.
\end{enumerate}
\mn
[Why?  Let $N_{i,j} = N_j \rest (M_j \cup \{s_j,t_j\})$.]

Easily {\bf tp}$(t_j,M_i,N_{i,j}) \le p_j$ and 
{\bf tp}$(t_j,M_i,N_{i,j}) = p_i$
by the claim \ref{b2.6} above.]
\mn
\begin{enumerate}
\item[$(*)_3$]   there is no $p \in \bold S^1(M_\theta)$ such that $i
< \theta \Rightarrow p \rest M_i = p_i$.
\end{enumerate}
\mn
Why?  We prove more: 
\mn
\begin{enumerate}
\item[$(*)_4$]   there is no $(N,t)$ such that
\begin{enumerate}
\item[$(a)$]  $M_\kappa \le_{\gk} N$
\sn
\item[$(b)$]  $t \in I^N$
\sn
\item[$(c)$] $(\forall \varphi \in \Gamma^{M_\kappa})[\varphi R^N_1 t]$.
\end{enumerate}
\end{enumerate}
\mn
[Why?  As then $\Gamma_\theta = \Gamma^M$ has a model contradiction to
an assumption.]  
\end{PROOF}

So e.g.
\begin{conclusion}
\label{b2.13}   
If $\theta > \kappa$ is regular with
no $\kappa^+$-complete uniform ultrafilter on $\theta$ and $\lambda =
2^\theta$, \then \, ${\gk}$ is not $(\lambda,\theta)$-sequence-compact.
\end{conclusion}

\begin{remark}
Recall if $D$ is an ultrafilter on $\theta$ then
min$\{\sigma':D$ is not $\sigma'$-complete$\}$ is $\aleph_0$ or a
measurable cardinality.
\end{remark}

\begin{PROOF}{\ref{b2.13}}
(Well known).

Let $M$ be the model with universe $2^\theta,P^M_0 = \theta$ and $R^M
\subseteq \theta \times \lambda$ be such that $\{\{\alpha < \lambda:\alpha
R^M \beta\}:\beta < \lambda\} = {\cP}(\theta),<^M$ the well
ordering of the ordinal on $\lambda$ the vocabulary has cardinality
$\kappa$ and has elimination of quantifiers and Skolem functions.

Let $\Gamma_i = \text{ Th}(M,\beta)_{\beta < \lambda} \cup\{\alpha <
c:\alpha < \theta\}$ ($c$ a new individual constant),
then $\langle \Gamma_i:i \le \theta\rangle$ is as\footnote{or directly
as $\Gamma_i$ has Skolem functions} required in \ref{b2.17} below
hence \ref{b2.11} apply.  
\end{PROOF}

\begin{conclusion}
\label{b2.17}   In Claim \ref{b2.11} if $\lambda =
\lambda^\kappa$ then we can allow $\langle \Gamma_i:i \le
\theta\rangle$ to be a sequence of theories in 
$\bbL_{\kappa^+,\kappa^+}(\tau),\tau$  any 
vocabulary of cardinality $\le \lambda$.
\end{conclusion}

\begin{PROOF}{\ref{b2.17}}
Without loss of generality we can add Skolem functions
(each with $\le \kappa$ places) in particular.  So $\Gamma_i$ becomes
universal and adding propositional variables for each quantifier free
sentence and writing down the obvious sentences, we get a set of
propositional sentences, we get $\Gamma_i$ as there.
\end{PROOF}

I think we forgot
\begin{observation}
\label{2b.21}  If $\lambda \ge \kappa \ge \theta =
\cf(\theta)$ \then \, the condition in \ref{b2.11} holds.
\end{observation}

\begin{proof}  Just let $\Gamma_0 = \{\bigvee_{i < \theta} \neg
p_i\},\Gamma_i = \Gamma_0 \cup \{p_j:j<i\}$.
\end{proof}

\begin{conclusion}
\label{2b.22}  
1) $\bold C_\kappa = \{\theta:\theta = \text{ cf}(\theta)$ 
and for every $\lambda$ and a.e.c. ${\gk}$ with LST$({\gk}) 
\le \kappa,|\tau_{\gk}| = \kappa$ have
$(\lambda,\theta)$-compactness of type$\}$ is the class
$\{\theta:\theta = \text{ cf}(\theta) > \kappa$ and there is a uniform
$\kappa^+$-complete ultrafilter on $\theta\}$.

\noindent
2) In $\bold C_\kappa$ we can replace ``every $\lambda$" by $\lambda =
   2^\theta + \kappa$.
\end{conclusion}

\begin{PROOF}{\ref{2b.21}}
Put together \ref{b2.13},\ref{b2.26}.
\end{PROOF}

Of course, a complimentary result (showing the main claim is best
possible) is:
\begin{claim}
\label{b2.23} 
If ${\gk}'$ is an a.c.c., {\rm LST}$({\gk}') \le \kappa$ 
and on $\theta$ there is a uniform
$\kappa^+$-complete ultrafilter on $\theta$ and $\theta$ is regular
and $\lambda$ any cardinality \then \, ${\gk}'$ has 
$(\lambda,\kappa)$-compactness of types.
\end{claim}

\begin{PROOF}{\ref{b2.23}}
Write down a set of sentences on $\bbL_{\kappa^+,\kappa^+}
(\tau^+_{\frak k})$ expressing the demands.

Let $\langle M_i:i \le \theta\rangle$ be $<_{\gk}$-increasing
continuous, $\|M_i\| \le \lambda,p_i = \text{\bf tp}_{\gk}
(a_i,M_i,N_i)$ so $M_i \le_{\gk} N_i$ such that $i<j < \theta
\Rightarrow p_i = p_j \rest M_i$.  Without loss of generality $\|N_i\|
\le \lambda$.

Let $\langle N_{i,j,\ell}:\ell \le n_{i,j,\ell}\rangle,\pi_{i,1}$
witness $p_i = p_j \rest M_i$ for $i < j < 0$ (i.e. $M_i \le_{\gk}
N_{i,j,\ell}$ (\wilog \, $\|N_{i,j,\ell}\| \le \lambda$), $N_{i,j,0} =
N_i,a_i \in N_{i,j,\ell},\bigwedge\limits_{\ell < n_{i,j,\ell}}
(N_{i,j,\ell} \le_{\gk} N_{i,j,\ell+1} \vee N_{i,j,\ell +1}
\le_{\gk} N_{i,j,\ell}$ and $\pi_{i,j}$ be an isomorphism from
$N_j$ onto $N_{i,j,n_{i,j}}$ over $M_i$ mapping $a_j$ to $a_i$.

Let $\tau^+ = \tau \cup\{F_{\varepsilon,n}:\varepsilon <
\kappa,n<\omega\}$, arity$(F_{\varepsilon,n}) = n$.  Let $\langle
M^+_i:i \le \theta\rangle$ be $\subseteq$-increasing, $M^+_i$ a
$\tau^+$-expansion of $M_i$ such that $u \subseteq M^+_i \Rightarrow
M_i \rest c \ell_{M^+_i}(u) \le_{\frak k} M_i$.  Similarly
$(N^{+,\varepsilon}_{i,j,\ell}:\ell \le n_{i,j,\ell}); \varepsilon =
1,\ell$ such that $N^{+,\varepsilon}_{i,j,\ell}$ is a $\tau^+$-expansion
of $N_{i,j,\varepsilon}$ as above such that $(\forall \ell <
n_{i,j,\ell})(\exists \varepsilon \in
\{1,2\})(N^{+,\varepsilon}_{i,j,\ell} \subset
N^{+,\varepsilon}_{i,j,\ell +1} \vee N^{+,\varepsilon}_{i,j,\ell +1}
\subseteq N^{+,\varepsilon}_{i,j,\ell})$.

Now write down a translation of the question, ``is there $p$ such that
..."
\end{PROOF}

\begin{claim}
\label{b2.26}  If $D$ is a uniform $\kappa$-complete
ultrafilter on $\theta,\langle M_i:i \le \theta \rangle$ is
$\le_{\gk}$-increasing continuous, $p_i \in \mathscr{S}^\alpha_{\gk}(M_i)$
as witnessed by $(N_i,a_i)$ for $i < \kappa,p_i = p_j \rest
M_i$ for $i<j< \kappa$ as witnessed by $(\pi_i,\langle
N_{i,j,\ell}:\ell \le m_{i,j}\rangle$ as in the proof above.

\noindent
1) There is $p_\kappa \in \bold S^\alpha(M_\theta)$ such that $i <
\theta \Rightarrow p_\kappa \rest M_i$.

\noindent
2) In fact for each $i < \kappa$ let ${\cU}_i \in D$ be such that
$i < j \in {\cU}_i \Rightarrow n_{i,j} = n^*_i$.  Let
  $N_{i,\kappa,\ell} = \prod\limits_{j \in {\cU}_i} N_{i,j,\ell}/D$.
So $\langle N_{i,\kappa,\ell}:\ell \le n^*_\ell\rangle$ are as 
above.   Let $M = \prod\limits_{i < \kappa} M_i/D,\pi_{i,\kappa} 
= \prod\limits_{j \in {\cU}_i} \pi_{i,j}/D$, etc.
\end{claim}
\newpage

\section {On some stability spectrums of an a.e.c.} 

\begin{convention}
\label{a1.5}  ${\gk}$ is an a.e.c. with amalgamation.
\end{convention}

\begin{definition}
\label{a1.7}  For $\theta \ge \text{ LST}({\gk})$.
We say ${\gk}$ is $(\lambda,\theta)$-stable \when \, $M \in 
K^{\gk}_\lambda \Rightarrow |\bold S(M)/E^\theta_M| \le \lambda$ where

\[
p E^\theta_M q \Leftrightarrow (\forall N)(N \le_{\gk} M \wedge
\|N\| \le \theta \Rightarrow p \rest N = q \rest \theta).
\]
\end{definition}

\begin{theorem}
\label{a1.11}  Fixing $\theta$ the class 
$\{\lambda:{\gk}$ is $(\lambda,\theta)$-stable$\}$; 
behave as in \cite{Sh:3}.
\end{theorem}

\begin{remark}
\label{a1.13}  See \cite{Sh:734} = \cite[V,\S7]{Sh:h} or
\cite{Sh:702} if not covered.
\end{remark}

\begin{definition}
\label{a1.17}    
$\kappa_\theta({\gk}) :=
\text{ Min}\{\kappa \le \theta^+$: there is no sequence $\langle M_i:i
\le \kappa\rangle$ which is $\le_{\gk}$-increasing continuous, $\|M_i\|
\le \theta$ and $p \in \mathscr{S}(M_\kappa)$ such that $p \rest M_{i+1}$
strongly $(\theta)$-split over $M_i\}$.
\end{definition}

\begin{claim}
\label{a1.19}  
1) If $\lambda > 2^\theta$ and ${\gk}$ is not
$(\lambda,\theta)$-stable \then \, for some $\kappa \le
\theta^+$ satisfying $\lambda^\kappa >\lambda$ we have
$\kappa < \kappa_\theta({\gk})$.

\noindent
2) If $\lambda > \theta,\lambda^\kappa > \lambda$ \then \, $\kappa <
\kappa_\theta({\gk})$ \then \, ${\frak k}$ not $(\lambda,\theta)$-stable.
\end{claim}

\begin{conclusion}
\label{a1.21}  $(-,\theta)$-stability spectrum - behave as in \cite{Sh:3}.
\end{conclusion}

\begin{discussion}
\label{a1.23}  
We can look at $\lambda \in [\theta,2^\theta)$ 
using splitting rather than strongly splitting.
\end{discussion}

It seems to me the main question is
\begin{question}
\label{a1.26}  Assume $(\exists \theta \ge 
\text{ LS}({\gk})(\kappa_\theta({\frak k}) > \aleph_0)$.

What can you say on Min$\{\theta:\kappa_\theta({\gh}) >
\aleph_0,\theta \ge \text{ LST}({\gk})\}$?
\end{question}

\begin{question}
\label{a1.31}  Assume GCH can we find an a.e.c. ${\gk}$ such that:
$(\forall \theta \ge \text{ LST}({\gk}))(\kappa_\theta({\gk})
= \aleph_0)$ \underline{but} unstable in every regular $\lambda > \text{
LST}({\gk})$?
\end{question}

\bibliographystyle{alphacolon}
\bibliography{lista,listb,listx,listf,liste,listy}

\end{document}